\documentclass{article}
\PassOptionsToPackage{square, numbers, compress}{natbib}
\usepackage{caption}
\usepackage{arxiv}
\usepackage{wrapfig}
\usepackage{float}
\usepackage{longtable}
\usepackage{natbib}
\usepackage{makecell}
\usepackage{subcaption}
\usepackage[dvipsnames]{xcolor}
\usepackage{tcolorbox}
\newtcolorbox{mycontentbox}[1][]{
    colback=gray!15,         
    colframe=black,          
    boxrule=1pt,           
    arc=1mm,                 
    left=6pt,
    right=6pt,
    top=6pt,
    bottom=6pt,
    #1
}

\usepackage{multirow}
\usepackage[utf8]{inputenc} 
\usepackage[T1]{fontenc}    
\usepackage{hyperref}       
\usepackage{url}            
\usepackage{booktabs}       
\usepackage{graphicx}      
\usepackage{amsfonts}       
\usepackage{nicefrac}       
\usepackage{microtype}      
\usepackage{xcolor}         
\usepackage{amsmath}
\usepackage{tabularx}
\usepackage{pifont}
\title{BenLOC: A Benchmark for Learning to Configure MIP Optimizers}
\usepackage{authblk}
\author{Hongpei Li$^{\dagger, }$}
\author{Ziyan He$^{\dagger, }$}
\author{Yufei Wang$^{\dagger, }$}
\author{Wenting Tu$^{*, }$}
\author[2]{Shanwen Pu}
\author[3]{Qi Deng}
\author[3]{Dongdong Ge}

\affil[1]{Shanghai University of Finance and Economics}
\affil[2]{Cardinal Operations}
\affil[3]{Shanghai Jiao Tong University}

\begin{document}
\maketitle
\footnotetext{$\dagger$ These authors contributed equally.}
\footnotetext{* Corresponding author.}
\footnotetext{$\ddagger$ Email: ishongpeili@gmail.com (Hongpei Li), heziyan@stu.sufe.edu.cn (Ziyan He), wyufff@163.sufe.edu.cn (Yufei Wang), tu.wenting@mail.shufe.edu.cn (Wenting Tu), 2019212802@live.sufe.edu.cn (Shanwen Pu), qdeng24@sjtu.edu.cn (Qi Deng), ddge@sjtu.edu.cn (Dongdong Ge).}

\begin{abstract}
The automatic configuration of Mixed-Integer Programming (MIP) optimizers has become increasingly critical as the large number of configurations can significantly affect solver performance. Yet the lack of standardized evaluation frameworks has led to data leakage and over‑optimistic claims, as prior studies often rely on homogeneous datasets and inconsistent experimental setups. To promote a fair evaluation process, we present BenLOC, a comprehensive benchmark and open-source toolkit, which not only
offers an end-to-end pipeline for learning instance-wise MIP optimizer configurations, but also standardizes dataset selection, train-test splits, feature engineering and baseline choice for unbiased and comprehensive evaluations of the effectiveness. Leveraging this framework, we conduct an empirical analysis on five well-established MIP datasets and compare classical machine learning models with handcrafted features against state-of-the-art deep-learning techniques.
The results demonstrate the importance of datasets, features and baseline criteria proposed by BenLOC and the effectiveness of BenLOC in providing unbiased and comprehensive evaluations.
\end{abstract}

\section{Introduction}

Mixed-Integer Programming (MIP) is a widely used mathematical optimization framework for solving real-world problems such as production scheduling, transportation planning, and revenue management~\citep{chen2011applied}. However, the NP-hard nature of MIP poses significant challenges, even for advanced MIP optimizers like GUROBI, CPLEX, and SCIP, especially when addressing large-scale, real-time industrial problems asking for rapid response~\citep{anand2017comparative}. Therefore, improving the efficiency of MIP solvers remains a popular topic of operations research.

Since MIP instances often share patterns and characteristics, machine learning (ML) approaches have emerged to enhance state-of-the-art solvers~\citep{khalil2022mip, gasse2019exact,HUANG2022108353}. While many studies use GNNs to optimize branch-and-bound and cutting plane methods, integrating ML into core algorithms is challenging due to the complexity of modifying well-optimized solvers. Thus, recent research has increasingly focused on ML-based instance-wise configuration, which operates externally without altering the solver’s internal routines~\citep{Iommazzo_2020}.

MIP solvers offer configurable settings for components like relaxations, heuristics, cutting planes, and branching strategies, all of which significantly impact solution quality and efficiency. Despite existing default settings and tuning tools, identifying optimal configurations for specific instances remains challenging and often requires extensive manual effort and expertise~\citep{inbook}.

Although there are some works focusing on the automatic configuration of optimizers, the lack of a unified evaluation standard and benchmark has limited their effectiveness in the following ways: (1) coverage of datasets is not comprehensive, as some rely on homogeneous datasets that fail to represent real-world variability; (2) data splitting methods, such as permutation-based splits, introduce data leakage, thereby inflating performance estimates; and (3) the lack of a consistent, data-driven baseline raises difficulties in fair comparisons.

With a thorough analysis of possible exaggerated performance results brought by existing methodologies, we propose a standard benchmark for ML-based optimizer automatic configuration by explicitly demonstrating: (1) a data-driven strategy for selecting datasets and baselines for unbiased comparisons; (2) the potential of underutilized dynamic handcrafted features, along with careful consideration of the extra time costs incurred when using features from solving process; also, we publish optimizer automatic configuration libraries and datasets built on our proposed benchmark, to promote future evaluations and comparisons in this field.

We summarize our primary contributions as follows. First, we present \textbf{the first comprehensive collection of datasets} tailored for learning MIP optimizer configurations, along with well-defined selection criteria to ensure their representativeness and benchmarking suitability. Second, we establish \textbf{a standardized experimental framework} that enables objective and consistent evaluation, mitigating exaggeration and bias. Third, we conduct \textbf{the most comprehensive benchmarking of competing methods to date}, from both feature and algorithmic perspectives, providing a clear comparison of their respective strengths and limitations. Finally, we release \textbf{an open-source toolset for optimizer configuration}, which includes datasets, feature extractors, and machine learning methods ranging from traditional Machine Learning to Deep Learning techniques.

\section{Literature Review}

Instance-Specific Algorithm Configuration (ISAC) \citep{ISAC} is a foundational approach using learning models to predict solver configurations, which clusters instances based on hand-crafted features (e.g., problem size, constraint types) using G-means clustering and applies a genetic algorithm (GGA) to find optimal configurations within each cluster. Hydra-MIP \citep{xu2011hydra} extends ISAC by iteratively selecting multiple configurations, enhancing performance on diverse MIP instances. Unlike these black-box approaches, \citet{Iommazzo_2020} formulate the configuration task as a mixed-integer nonlinear program (MINLP) using Support Vector Regression (SVR) to learn a performance map, achieving notable improvements in Hydro Unit Commitment (HUC) problems.

\begin{wraptable}{r}{0.6\linewidth}   
    \vspace{0pt}                       
    \centering
    \caption{\scriptsize Literature Review. Dataset abbreviations: SC = Set Cover, CA = Combinational Auction, CFL = Capacitated Facility Location, IS = Independent Set, PAC = Packing, BIN = Binary Packing, PLA = Planning, MC = Maximum Cut, MRP = Multiple RGV Path, COR = CORLAT, MIK = MIK, IP = Item Placement, LB = Load Balancing, ANO = Anonymous, MIRP = Maritime IRP, CLU = CLUREG, CLUUR = CLUREGURCW, ISAC = ISAC (new), HUC = Hydro Unit Commitment, CMIP = Common MIP, MIP = MIPLIB, PP = Production Planning, GMI = General MIP.}

    \label{literature}
 {
  \fontsize{7}{8}\selectfont
    \setlength\tabcolsep{0.5pt}
\begin{tabular}{lcccccccc}
    \toprule
    \multirow{2}{*}{\textbf{Liter.}} & \multicolumn{2}{c}{\textbf{Method}} & \multicolumn{2}{c}{\textbf{Dataset}} & \multirow{2}{*}{\textbf{Baseline}} & \multicolumn{2}{c}{\textbf{Feature}}  & \multirow{2}{*}{\textbf{\makecell{Imp \\ Default}}} \\

    \cmidrule(lr){2-3} \cmidrule(lr){4-5} \cmidrule(lr){7-8}
    & \textbf{ML} & \textbf{DL} & \textbf{Homo.} & \textbf{Hetero.} &   & \textbf{Handcraft} & \textbf{Graph} &    \\
    \midrule
    \citep{gasse2019exact} & & \makecell{GNN} & \makecell{SC\\CA, CFL \\ IS} & \ding{55} & \makecell{default, FSB \\SVMRank \\ LMART, Trees} & & \makecell{GNN} & \makecell{-74\%-36\%} \\
    \midrule
    \citep{tang2020reinforcement} & \makecell{RL} & & \makecell{PAC, BIN \\ PLA, MC} & \ding{55} & \makecell{Heuristics:\\ Random \\ MV, MNV,LE} & \makecell{static} &  & \\
    \midrule
    \citep{l2p} & & \makecell{GNN} & \makecell{SC, IS, MRP \\ COR, MIK \\ IP, LB, ANO, MIRP} & \ding{55} & \makecell{default\\Random\\ SMAC3\\FBAS} & & \makecell{GNN}  & \makecell{0-50\%} \\
    \midrule
    \citep{xu2011hydra} & \makecell{Tree} & & \makecell{CLU, CLUUR\\ ISAC} & MIX & \makecell{default,\\PARAMILS,\\MIPzilla} & \makecell{static\\dynamic} & \makecell{graph}  & \makecell{0-18\%} \\
    \midrule
    \citep{Iommazzo_2020} & \makecell{SVR} & & \makecell{HUC} & \ding{55} & \makecell{default} & \makecell{static} &  & \\
    \midrule
    \citep{valentin2022instancewisealgorithmconfigurationgraph} & & \makecell{GNN} & \makecell{IP, ANO} & \ding{55} & \makecell{default} & & \makecell{GNN}  & \\
    \midrule
    \citep{clustering} & & \makecell{GNN} & \makecell{CA, CFL, SC, IS} & \makecell{CMIP, MIP} & \makecell{default\\SMAC\\ ISAC} & \makecell{static} & \makecell{GNN}  & \makecell{-63\%-22.9\% \\ (MIPLIB)} \\
    \midrule
    \citep{HUANG2022108353} & & \makecell{MLP} & \makecell{SC, KP, PLA \\ GMI, PP} & \ding{55} & \makecell{Manually \\designed \\ heuristics} & \makecell{static\\dynamic} &  & \makecell{28\%-52\%} \\
    \midrule
    \citep{berthold2022learning} & \makecell{Rfr} & & \makecell{FICO} & \makecell{MIP} & \makecell{Always LC} & \makecell{static\\dynamic} &  & \makecell{9-11\%} \\
    \bottomrule
\end{tabular}
}
\vspace{0pt}
\end{wraptable}

Recent studies also leverage graph-based representations to better capture structural information in MIP instances. \citet{valentin2022instancewisealgorithmconfigurationgraph} propose a GNN-MLP hybrid model that represents MIP problems as bipartite graphs, enabling end-to-end prediction of optimal SCIP configurations. Similarly, \citet{clustering} introduce DGCAC, which uses random-walk-based graph embeddings and auto-encoders to cluster MIP instances and apply SMAC \citep{10.1007/978-3-642-25566-3_40} for cluster-wise configuration, outperforming manual feature-based methods. For presolving configurations, \citet{l2p} develop a framework, named L2P, that combines neural networks with simulated annealing, using residual networks to capture presolver interactions, resulting in significant performance gains on industrial benchmarks.
Table \ref{literature} summarizes recent advancements in learning-based methods for MIP configuration. Despite their potential, most studies have been tested on simple synthetic or homogeneous datasets, with limited success on more complex benchmarks like MIPLIB\citep{l2p}. Moreover, performance comparisons often focus only on default solver settings, overlooking more rigorous comparisons against optimal configurations for specific instances or datasets. This lack of comprehensive evaluation raises concerns about the generalizability of these methods across different solvers and datasets.

To address these gaps, this paper introduces a standardized benchmarking framework that includes diverse datasets and well-defined evaluation metrics, providing a more accurate assessment of ML-based configuration methods in MIP solvers.

\section{Preliminary: Algorithm Configuration for MIP Optimizers}

In the context of MIP optimizers, Algorithm Configuration (AC) aims to select a configuration based on solver performance metrics, such as runtime, primal gap, or dual gap. Currently, algorithms implementing Algorithm Configuration (AC) can be classified into two types based on the granularity of the obtained optimal configuration: Per-Instance (PI) and Per-Dataset (PD).  The PI approach seeks to determine the optimal configuration $\theta(x) \in \Theta$ for each instance $x \in \mathcal{X}$, requiring prior known instance-specific features. In contrast, the PD approach identifies a single configuration $\theta^*$ that performs well on average across the entire dataset, without considering individual instance variability \cite{Iommazzo_2020}. While the PD approach simplifies the configuration process, it may struggle when algorithm performance varies significantly across instances.

Due to the considerable variation in solver performance across instances belonging to the same problem class, most studies focus on the PI approach. Specifically, they leverage historical instances \( \{x_1, \ldots, x_n\} \) and their corresponding optimal configurations \( \{\theta_1^*, \ldots, \theta_n^*\} \) to train a predictive model \( \mathcal{M} \). The objective is to identify the instance-specific optimal configuration \( \theta^*(x) \) that minimizes the performance measure: \(
\theta_{\text{pred}} := \arg\min_{\theta \in \Theta} \mathcal{M}(x)\), where \( \mathcal{M}(x) \) is the predicted configuration for a new instance \( x \).

Conversely, the PD approach, which identifies a single configuration $\theta^*$ that performs well on average across the MIP dataset, is often overlooked. In the following section, we consider both the PI best configuration and PD best configuration, emphasizing their respective strengths and limitations.

\section{BenLOC: A Benchmark of Learning to MIP Optimizer Configuration}
We present BenLOC which outlines criteria for selecting datasets with potential for improvement, proposes data splitting strategies to prevent leakage during data augumentation, and defines robust baselines for fair evaluation. Together, these components provide unbiased assessments and a strong foundation for advancing ML-based MIP optimizer configuration.

\subsection{Learning-Based MIP Optimizer Configuration Pipeline}

\begin{figure*}[htbp]
    \centering
    \includegraphics[width=1.0\linewidth]{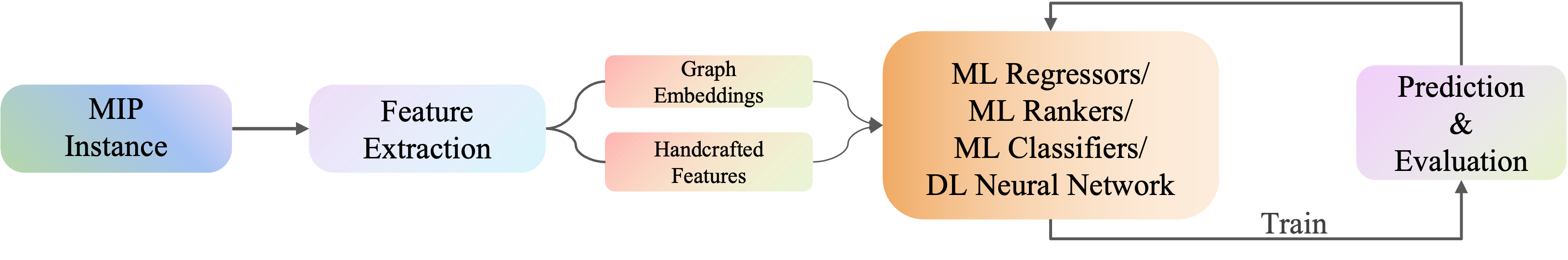}
    \caption{MIP Optimizer Configuration Workflow}
    \label{fig:workflow}
\end{figure*}

When using learning techniques for MIP optimizer configuration, the goal is to train a model to predict the best configuration for each new instance. First, relevant features are extracted from MIP instances, which can be either handcrafted features derived from basic instance properties or graph embeddings that capture structural information. These features are then fed into a model, which can be an ML regressor, a ranker, or a deep learning network. The objective is to predict configurations that minimize certain evaluation metrics, with BenLOC specifically focusing on solving time. The overall workflow is illustrated in Figure~\ref{fig:workflow}.

\subsection{Dataset Collection and Selection Criteria}
\par\textbf{Homogeneous or Heterogeneous?} Datasets used to evaluate MIP solvers can be categorized into homogeneous and heterogeneous types. Homogeneous datasets, like Set Covering and Load Balancing~\cite{gasse2022machine}, focus on a single MIP problem structure, while heterogeneous datasets like MIPLIB~\cite{MIPLIB2017} contain diverse instances with varying constraints, scales, and objectives. However, as shown in Table~\ref{literature}, most works on learning the MIP optimizer configuration focus on homogeneous datasets, while other studies considering MIPLIB often struggle to outperform default configurations. Thus, testing on heterogeneous datasets can better show the potential of proposed models. Moreover, the instances met by solvers in practical applications typically stem from various industry backgrounds with different structures, making them more similar to heterogeneous datasets. Therefore, research on ML-based optimizer configuration should include heterogeneous datasets such as MIPLIB.

\par\textbf{The Value of Using ML and DL-based Algorithms} ML and DL-based algorithms aim to learn a model that maps instance features to Per-Instance (PI) configurations, thereby pursuing optimal PI configurations.  To validate the effectiveness of a learning-based algorithm, the selected datasets should be those for which the optimal Per-Instance (PI) configuration yields significantly better performance than the optimal Per-Dataset (PD) configuration. Otherwise, it becomes unclear whether the algorithm fails to accurately predict the optimal PI configuration or whether the optimal PI configuration itself offers no substantial advantage.
\begin{wrapfigure}{r}{0.43\linewidth}   
    \vspace{-6pt}                       
    \centering
    \captionof{table}{Improvement on Different Datasets}
    \label{ImpDataset}
    {\fontsize{6}{6}\selectfont
    \setlength\tabcolsep{2pt}
    \begin{tabularx}{\linewidth}{lcccc}  
        \toprule
        \textbf{Dataset} & \textbf{Instance Cnt.} & \textbf{PD Best (\%)} & \textbf{PI Best (\%)} & \textbf{Imp. UB (\%)} \\
        \midrule
        Set Cover & 3994 & 25.67\% & 26.53\% & \textbf{0.86\%} \\
        Independent Set & 1604 & 1.84\% & 9.45\% & 7.61\% \\
        NN Verification & 3104 & 6.12\% & 37.84\% & 31.72\% \\
        Load Balance & 2286 & 5.59\% & 25.08\% & 19.49\% \\
        MIPLIB 2017 & 1065 & 0.67\% & 34.66\% & 33.99\%\\
        \bottomrule
    \end{tabularx}
    }
    \parbox{\linewidth}{\tiny\textit{Data is derived from the solving results of Cardinal Optimizer (COPT \citep{copt} on 2 configurations for MIPLIB and 4 configurations for other datasets}}
    \vspace{-6pt}                       
\end{wrapfigure}
Table \ref{ImpDataset} illustrates the improvement of the Per-Instance best over the Per-Dataset best across five datasets. Note that the PI-best serves as the ground-truth, representing the upper bound of performance that learning-based algorithms can potentially achieve. Notably, while the commonly used Set Cover dataset in MIP optimizer configuration shows substantial improvement over the default configuration, the improvement relative to the Per-Dataset best configuration is minimal. This suggests that the Per-Dataset best configuration is already highly effective, leaving limited room for ML and DL-based configuration tuning to provide further gains. This observation is further supported by the experimental results in Table \ref{exp2}.


\subsection{Data Splitting on Augmented Datasets for Unexaggerated Evaluation}

Given the limited number of training instances in MIPLIB, \cite{berthold2022learning, Berthold_Hendel_2021} propose a permutation-based data augmentation strategy that generates additional training data by creating multiple permutations of existing instances through row and column swaps in the constraint matrix. This method is valid since while the problem stays unchanged mathematically, it is well-known that such problem permutations can have a dramatic impact on the performance.

When performing a train-test split on a dataset augmented using permutation techniques, it is important to ensure that the splitting criterion prevents any instance in the test set from having its permutations that also appear in the training set. The reason is that permutations are inherently related. If permutations of a test sample appears in the training set, they can significantly influence the prediction for that test sample, and this influence is often positive. However, in practical applications, when predicting a test sample, we do not generate its permutations and obtain ground-truth labels for them beforehand. In other words, in real-world scenarios, test samples do not have their permutations included in the training set. Therefore, BenLOC strongly recommends a split-by-instance strategy, where each unique problem is fully allocated to either the training or test set when using permutations for data augmentation, and does not advise using split-by-permutation method, is likely to cause data leakage.  In Section \ref{finding:1}, we will provide experimental results showing that improper data splitting can result in an inflated assessment of the algorithm performance.


\subsection{Features Engineering}\label{feature}
In this section, we present the feature engineering process in BenLOC, covering both handcraft features and graph-based features extracted by GNN. In particular, we will discuss the importance of dynamic features for the predictability of the optimal configuration and propose guidelines for their use.

\subsubsection{Handcraft Features}

Handcrafted features in MIP are domain-specific attributes extracted from the problem file and the optimizer log, which can be divided into static features and dynamic features. Static features capture the mathematical structure, such as the number of constraints and constraint types, and remain unchanged throughout the solving process. Dynamic features, however, evolve as the solving process progresses, as shown in Figure \ref{fig:copt}. For example, the dynamic feature \# Sepa only gains meaningful values after the initial root LP is completed, remaining zero up until that point.

\textbf{Importance of Dynamic Features }However, static features like problem size and constraint types provide limited insight into instance difficulty, as even minor changes, such as matrix permutations, can significantly affect solver performance without altering these features. Thus, the introduction of dynamic features is necessary. In BenLOC, we propose three different feature settings: (1) only static features, (2) static features combined with dynamic features up to the first root LP, and (3) dynamic features up to the completion of all root nodes. Their effectiveness is demonstrated in the experimental section.

\begin{minipage}{0.69\textwidth}
\centering
\captionof{table}{Static and Dynamic Features}
\resizebox{\linewidth}{!}{
\begin{tabular}{|c|l|c|l|}
\hline
\multicolumn{2}{|c|}{\textbf{Static Feature}} & \multicolumn{2}{c|}{\textbf{Dynamic Feature}} \\ \hline
\textbf{Category} & \textbf{Feature} & \textbf{Category} & \textbf{Feature} \\ \hline
\multirow{3}{*}{Matrix} & Rows  & \multirow{3}{*}{Presolving} & PresolRows  \\ \cline{2-2} \cline{4-4}
& Columns  &  & PresolColumns  \\ \cline{2-2} \cline{4-4}
& NonZeros  &  & PresolIntegers  \\ \cline{1-2} \cline{3-4}
\multirow{2}{*}{Variables} & Binaries  & \multirow{4}{*}{\shortstack{Global\\Cutting}}  & DualInitialGap  \\ \cline{2-2} \cline{4-4}
& Integers  &  & PrimalDualGap  \\ \cline{1-2} \cline{4-4}
\multirow{13}{*}{Constraints} & LessThan  &  & PrimalInitialGap \\ \cline{2-2} \cline{4-4}
& GreaterThan &  & GapClosed  \\ \cline{2-2} \cline{3-4}
& Equality  & \multirow{7}{*}{\shortstack{Finish First\\Root LP}}  & Active  \\ \cline{2-2} \cline{4-4}
& SetPartitioning &  & IntInf  \\ \cline{2-2} \cline{4-4}
& SetPacking  &  & GlbRed  \\ \cline{2-2} \cline{4-4}
& SetCovering &  & Gap  \\ \cline{2-2} \cline{4-4}
& Cardinality &  & Time  \\ \cline{2-2} \cline{4-4}
& KnapsackEquality &  & objective\_density \\ \cline{2-2} \cline{4-4}
& Knapsack  &  & Symmetries \\ \cline{2-2} \cline{3-4}
& KnapsackInteger & \multirow{7}{*}{\shortstack{Finish \\Root Node}}  & Nodes  \\ \cline{2-2} \cline{4-4}
& BinaryPacking &  & LPit/n  \\ \cline{2-2} \cline{4-4}
& VariableLowerBound &  & GlbFix  \\ \cline{2-2} \cline{4-4}
& VariableUpperBound &  & \#Cuts  \\ \cline{2-2} \cline{4-4}
& MixedBinary &  & \#MCP  \\ \cline{1-2} \cline{4-4}
\multirow{3}{*}{Scaling}  & Coefficient\_oom &  & \#Sepa  \\ \cline{2-2} \cline{4-4}
& RightHandSide\_oom &  & \#Conf  \\ \cline{2-2} \cline{3-4}
& Objective\_oom & & \\ \hline
\end{tabular}}

\label{Table Feature}
\end{minipage}
\hfill
\begin{minipage}{0.29\textwidth}
\centering
\includegraphics[width=\linewidth]{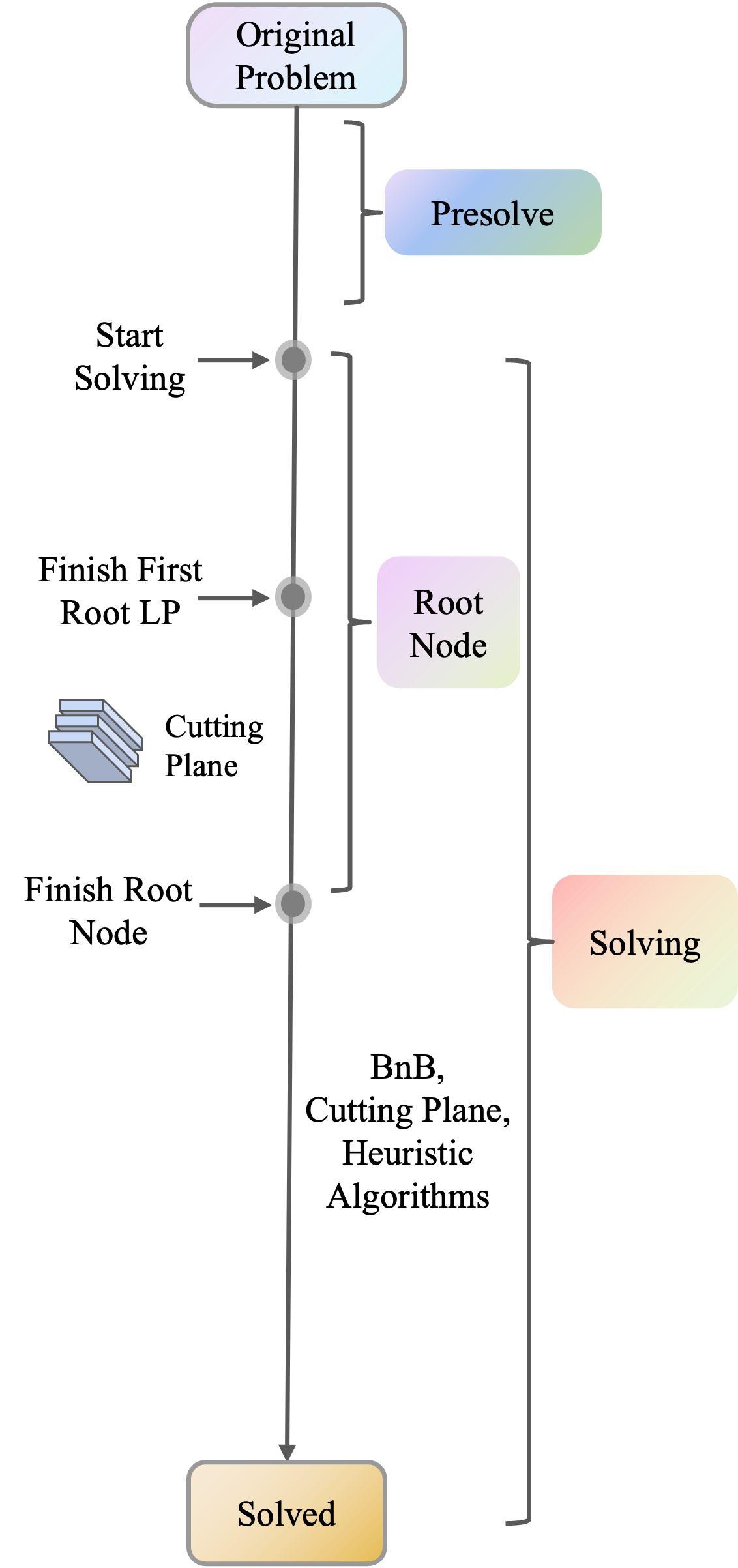}
\vspace{0.5em}

\captionof{figure}{COPT Solving Process}
\label{fig:copt}
\end{minipage}

\textbf{Extra Cost Overlooked} As the solver progresses, features obtained in later stages may provide stronger predictive power, but they can also introduce extra time costs when influenced by testing configurations. As demonstrated in Figure~\ref{fig:copt}, the MIP-solving process involves multiple distinct phases with various algorithms, each of which can be influenced by different configurations, altering the problem state. Therefore, when training with ML, an ideal approach is to utilize features that can be obtained before the configurations are applied to the MIP. For example, in the case of Treecut, which operates after completing the root node, using features extracted prior to this stage incurs no additional cost. In contrast, if the current testing configuration influences the root node solving process, such as a Subheuristic, then the configuration must be assigned before completing the root nodes. This implies that after prediction, optimizers need to revert to the first root LP point to assign the predicted configuration, resulting in an additional time cost to resolve the root node.  Ignoring these costs can lead to artificially inflated test performance (as shown in the experimental section).

\subsubsection{Graph-based Features}
Previous works show that a bipartite graph formulation can effectively represent an MIP instance \citep{gasse2019exact, knn-gnn}, where one set of nodes corresponds to variables and the other to constraints, and edges represent variable participation in constraints. This structure is well-suited for applying Bipartite Graph Neural Networks (GNNs), which can learn meaningful representations of variables and constraints by propagating and aggregating information over the graph. These learned embeddings serve as features that capture the structural properties and interactions within the MIP instance. We provide further architectural and implementation details in the Appendix \ref{gnn-feature}.

\subsection{Baselines in BenLOC}
Previous studies typically use the default configuration as the baseline, as it represents the optimizer’s general-purpose tuning \cite{l2p}. However, default configurations are designed for general-purpose use and lack the specificity required to address the unique characteristics of individual datasets and these generic settings often underperform compared to configurations tailored to the dataset's problem distribution.
Therefore, to achieve a more comprehensive and objective evaluation, a greater number of baseline methods for comparison are required. In BenLOC, we suggest considering the Per-Dataset best configuration as a stronger baseline, which corresponds to the configuration that performs best on average across a particular dataset. This is because the Per-Dataset best helps determine whether learning-based algorithms can effectively model the relationship between instances and configurations. Furthermore, the Per-Dataset best is a simple approach that can be applied in a wide range of scenarios, making it a representative baseline. Additionally, BenLOC provides comprehensive and generalizable features as well as implementations of ML and DL algorithms for learning-based configuration tuning, which can also serve as baselines for future research.

\section{Experiments and Empirical Findings}

\subsection{Experimental Setup}
\textbf{Datasets and Features} To comprehensively evaluate and obtain generalizable results, we utilize five datasets: four homogeneous datasets (Set Covering, Independent Set, NN Verification, Load Balancing) and one heterogeneous dataset (MIPLIB2017). Each instance is tested on COPT using 25 configurations, with the top three configurations based on average performance selected alongside the default configuration. For MIPLIB2017, we augment the dataset by testing on 10 different permutations, providing additional data diversity. Table \ref{ImpDataset} summarizes the dataset statistics. Feature extraction includes both handcrafted features (with dynamic features derived from COPT solving logs) and graph features, as outlined in Section \ref{feature}. All instances, along with their corresponding solving times and extracted features, are integrated into the BenLOC toolkit to enable replication and further analysis. Notably, while the experimental framework is based on COPT, it can be readily extended to other solvers, enhancing its applicability across various optimization platforms.

\textbf{Models} In our experiments, we implement both traditional machine learning techniques and deep learning models. The machine learning methods include regression models (Random Forest, LightGBM, XGBoost, GBDT), ranking models (LightGBM, XGBoost), and classification models (Random Forest, LightGBM, XGBoost, GBDT, KNN, SVM). Bayesian optimization is employed for ML model configuration. For deep learning methods, we employed TableNet, GNN-MLP \citep{valentin2022instancewisealgorithmconfigurationgraph} and L2P \citep{l2p} as representatives of general-purpose and MIP-specific deep learning techniques, respectively. See the appendix for details. For simplicity, only a subset of these methods is presented in Table \ref{exp2}. More comprehensive experimental results are presented in Table~\ref{exp2-all} in the Appendix.

\textbf{Evaluation Metrics} We compare the model results with both COPT's default configuration and the PD best configuration for each dataset. Improvement is calculated based on the shift geometric average time as follows: \(
\text{Improvement} = \frac{\text{Baseline} - \text{Preddict Time}}{\text{Baseline}}
\), where \text{Baseline} can be either the default configuration or the PD best configuration.




\subsection{Influence of Data-splitting Strategies }\label{finding:1}

We conduct experiments on the permutated MIPLIB dataset to evaluate the impact of permutation-based splitting and instance-based splitting. To ensure the robustness of our results and mitigate the impact of randomness, we perform 10 independent train-test splits and report the average performance.

\begin{minipage}{0.70\textwidth}
    \centering
    \captionof{table}{Performance Comparison of Different Data Splitting Methods}
    \resizebox{\linewidth}{!}{
    \begin{tabular}{lcccccccc} 
        \toprule
        & \multicolumn{4}{c}{\textbf{Split by Instance}} & \multicolumn{4}{c}{\textbf{Split by Permutation}} \\
        \cmidrule(lr){2-5} \cmidrule(lr){6-9} 
        \textbf{} & \textbf{Pred.} & \textbf{PD Best} & \textbf{PI Best} & \textbf{Imp.} & \textbf{Pred.} & \textbf{PD Best} & \textbf{PI Best} & \textbf{Imp.} \\
        \midrule
        Random Forest & 45.39 & 46.55 & 41.76 & 2.49\% & 43.59 & 45.98 & 36.34 & 4.34\% \\
        LightGBM  & 46.41 & 46.55 & 41.76 & 2.69\% & 43.56 & 45.98 & 36.34 & 4.41\% \\
        XGBoost   & 45.32 & 46.55 & 41.76 & 0.41\% & 43.63 & 45.98 & 36.34 & 4.26\% \\
        GBDT    & 45.96 & 46.55 & 41.76 & 1.27\% & 43.69 & 45.98 & 36.34 & 4.13\% \\
        Rfr Clf   & 46.05 & 46.55 & 41.76 & 0.96\% & 44.08 & 45.98 & 36.34 & 3.27\% \\
        LGBM Clf  & 46.08 & 46.55 & 41.76 & 0.90\% & 44.15 & 45.98 & 36.34 & 3.12\% \\
        XGB Clf   & 46.20 & 46.55 & 41.76 & 0.75\% & 43.97 & 45.98 & 36.34 & 3.51\% \\
        GBDT Clf  & 46.10 & 46.55 & 41.76 & 0.97\% & 44.34 & 45.98 & 36.34 & 2.70\% \\
        KNN Clf   & 46.15 & 46.55 & 41.76 & 0.86\% & 41.49 & 45.98 & 36.34 & 9.58\% \\
        SVM     & 46.24 & 46.55 & 41.76 & 0.67\% & 41.73 & 45.98 & 36.34 & 9.18\% \\
        \bottomrule
    \end{tabular}}
    \label{exp_split}

\end{minipage}
\hfill
\begin{minipage}{0.27\textwidth}
    \centering
        \includegraphics[width=\textwidth]{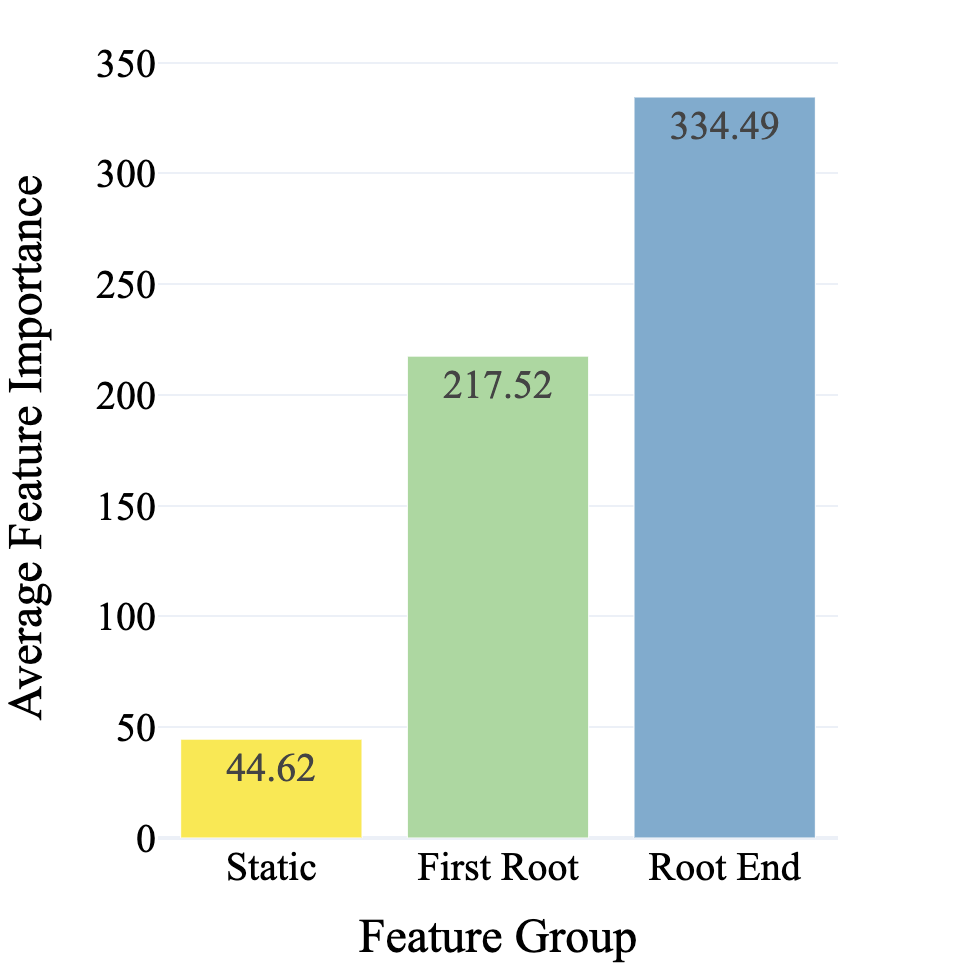}
        \captionof{figure}{Feature Importance}
        \label{fig:importance}
\end{minipage}

\begin{mycontentbox}
\textbf{Split-by-Permutation Leads to Inflated Performance} Permutation-based splits introduce data leakage by allowing similar problem structures in both training and testing sets, resulting in overestimated performance evaluations. Conversely, instance-based splitting with permutation-augmented data provides a more rigorous evaluation of ML-based optimizer configuration.
\end{mycontentbox}
    Results shown in Table~\ref{exp_split} confirm that permutation-based data splitting significantly inflates model performance, as shown. When using permutation-based splits, the test performance improves by up to 9.18\% by using KNN as the classifier compared with PD-best configuration. However, under instance-based splitting, the improvement drops to 2.69\% suggesting that permutation-based splits will actually introduce data leakage and over-optimistic performance evaluations.
    
\subsection{Influence of Dataset, Feature and Baseline Choice }\label{finding:2}

To thoroughly demonstrate the effectiveness of the standards proposed by BenLOC in terms of datasets, features, and baselines, we conducted experiments on multiple datasets, incorporating static, dynamic, and graph-based features, as well as a variety of ML and DL approaches. The results are presented in Table~\ref{exp2} and Table~\ref{exp_gnn}.
\begin{table}[htbp]
    \centering
    \caption{Performance Comparison with Different Datasets, Features and Baselines with Handcrafted Features}
    \label{exp2}
  \resizebox{\linewidth}{!}{
    \begin{tabular}{lcccccccccccc}
        \toprule
    &  \textbf{Predictor} & \textbf{Default} & \textbf{PD Best}
    & \multicolumn{3}{c}{\textbf{Static Feature}}& \multicolumn{3}{c}{\textbf{Static \& First Root Dynamic Feat.}} & \multicolumn{3}{c}{\textbf{Static \& Root End Dynamic Feat.}} \\
    \cmidrule(lr){5-7} \cmidrule(lr){8-10} \cmidrule(lr){11-13}
    & & &  & \textbf{Predict} & \textbf{Imp.}& \textbf{Imp.} & \textbf{Predict} & \textbf{Imp.}& \textbf{Imp.} & \textbf{Predict} & \textbf{Imp.}& \textbf{Imp.} \\
    & & &  & \textbf{Time}& \textbf{Default} & \textbf{PD Best} & \textbf{Time} & \textbf{Default} & \textbf{PD Best} & \textbf{Time} & \textbf{Default} & \textbf{PD Best} \\

    \bottomrule
        \midrule
        \multirow{5}{*}{\shortstack{Independent\\Set}}
        &Random Forest & 5.92 & 5.79 &5.73& 3.21\% & 0.93\% & 5.73& 3.21\% & 0.94\% & 5.73& 3.20\% & 0.96\% \\

         & LightGBM
          & 5.92 & 5.79 &5.75&2.79\% & 0.64\% & 5.74& 2.96\% & 0.81\% & 5.73& 3.13\% & \textbf{0.98\%} \\

        & XGB Ranker
          & 5.92 & 5.79 & 5.86 & 1.01\% & -1.21\% & 5.88 &0.55\% & -1.65\% & 5.86&1.01\% & -1.21\% \\

         & KNN Clf
          & 5.92 & 5.79 & 5.79 & 2.12\% & -0.05\% & 5.77 &2.46\% & 0.29\% & 5.75&2.79\%& 0.64\% \\

        &TableNet & 5.92 & 5.79 &5.75& 2.86\% & 0.72\% & 5.74& 2.93\% & 0.79\% & 5.73& 3.16\% & 0.93\% \\

        \midrule
        \multirow{5}{*}{\shortstack{Load\\Balance}}
         &Random Forest & 1123.94 & 1068.90 & 1057.40& 5.92\% & 1.08\% & 1051.60& 6.44\% & 1.62\% & 1034.41& 8.23\% & 3.23\% \\

        & LightGBM  & 1123.94 & 1068.90  & 1059.55& 5.73\% & 0.88\% & 1057.45& 5.92\% & 1.07\% & 1030.25& 8.34\% & \textbf{3.62\%} \\

        & XGB Ranker  & 1123.94 & 1068.90 & 1084.89& 3.47\% & -1.50\% & 1085.24& 3.44\% & -1.53\% & 1086.81& 3.30\% & -1.68\% \\

           & KNN Clf & 1123.94 & 1068.90 & 1070.50& 4.75\% & -0.15\% & 1069.35& 4.86\% & -0.05\% & 1035.90& 7.83\% & 3.09\% \\

       &TableNet  & 1123.94 & 1068.90 & 1064.70& 5.27\% & 0.39\% & 1048.99& 6.68\% & 1.86\% & 1035.90& 7.84\% & 3.09\% \\

        \midrule

        \multirow{5}{*}{\shortstack{NN\\Verification}}
          &Random Forest  & 6.98 & 6.50 & 6.59& 5.59\% & -1.58\% & 6.58& 5.73\% & -1.40\% & 6.26& 10.32\% & 3.69\% \\

        & LightGBM  & 6.98 & 6.50 & 6.62& 5.13\% & -1.89\% & 6.63& 4.99\% & -2.05\% & 6.25& 10.43\% & \textbf{3.80\%} \\

        & XGB Ranker  & 6.98 & 6.50 & 6.71& 3.84\% & -3.28\% & 6.74& 3.41\% & -3.74\% & 6.70& 3.98\% & -3.12\% \\

           & KNN Clf & 6.98 & 6.50 & 6.65& 4.70\% & -2.36\% & 6.58& 5.70\% & -1.28\% & 6.38& 8.57\% & 1.80\% \\

        &TableNet  & 6.98 & 6.50 & 6.62& 5.07\% & -1.95\% & 6.59& 5.59\% & -1.58\% & 6.28& 10.12\% & 3.39\% \\

         \midrule

        \multirow{5}{*}{Set Cover}
        &Random Forest & 2.89 & 2.14 &2.14&26.19\% & -0.03\% & 2.14& 26.20\% & -0.02\% & 2.13& 26.30\% & \textbf{0.04\%} \\

        & LightGBM
          & 2.89 & 2.14  &2.14&26.19\% & -0.03\% & 2.14& 26.20\% & -0.02\% & 2.15& 26.05\% & -0.12\% \\

        & XGB Ranker
          & 2.89 & 2.14 & 2.59 & 10.68\% & -21.06\% & 2.58 &10.75\% & -20.96\% & 2.58 &10.61\% & -21.17\% \\

          & KNN Clf
          & 2.89 & 2.14 & 2.14 & 26.02\% & -0.29\% & 2.14 & 26.04\% & -0.26\% & 2.14 & 26.19\% & -0.05\% \\

          &TableNet & 2.89 & 2.14 &2.14 &26.18\% & -0.07\% & 2.14& 26.19\% & -0.06\% & 2.13& 26.28\% & 0.02\% \\

           \midrule
        \multirow{5}{*}{MIPLIB}
         &Random Forest  & 46.55 & 46.55 & 46.62 & -0.15\% & -0.15\% & 46.50 & 0.12\% & 0.12\% & 45.39 & 2.49\% & 2.49\%\\

         &LightGBM  & 46.55 & 46.55 & 47.55 & -2.31\% & -2.31\%  & 46.41 & 0.31\% & 0.31\% & 45.32 & 2.69\% & \textbf{2.69\%}\\

         &XGB Ranker  & 46.55 & 46.55 & 47.64 & -2.34\% & -2.34\% & 47.23 & -1.46\% & -1.46\% & 46.97 & -0.90\% & -0.90\%\\

 & KNN Clf & 46.55 & 46.55 & 47.83 & -2.74\% & -2.74\% & 46.50 & 0.11\% & 0.11\% & 46.15 & 0.86\% & 0.86\% \\

    & TableNet  & 46.55 & 46.55 & 46.61 & -0.13\% & -0.13\% & 46.52 & 0.06\% & 0.06\% & 45.42 & 2.43\% & 2.43\%\\

\bottomrule
\end{tabular}
}
\end{table}

\begin{mycontentbox}
\textbf{Datasets Unsuitable for Learning-based Configuration Tuning} Some datasets inherently lack the improvement space between the Per-Dataset best configuration and the optimal configuration, leading to limited benefits from learning-based methods even when multiple models are applied. This highlights the importance of dataset choice in realizing performance gains through automated configuration.
\end{mycontentbox}

Recall that, as shown in Table~\ref{ImpDataset}, Set Covering exhibits a minimal performance improvement between the Per-Dataset best configuration and the optimal configuration, indicating limited potential for ML-based configuration tuning. Specifically, the improvement over the Per-Dataset best setting is only 0.86\%, which is significantly lower compared to heterogeneous datasets like MIPLIB (33.99\%) and NN Verification (31.72\%). Thus, as indicated in Table~\ref{exp2}, while ML-based approaches can reduce the gap to the Per-Dataset best setting, the improvement for Set Covering remains close to 0\% and does not become positive, in contrast to other datasets where ML-based tuning yields noticeable and consistently positive gains.

\begin{mycontentbox}
    \textbf{Dynamic Features Exhibit Significantly Stronger Predictive Power} Dynamic features play a crucial role in enhancing model performance, with deeper-stage features providing richer solver behavior signals and achieving superior predictive accuracy. Models trained with only static features show limited or negative improvements, while incorporating dynamic features, particularly those extracted at later stages (e.g., after the root node), significantly boosts predictive performance.
\end{mycontentbox}

    We evaluate the impact of dynamic features on model performance by comparing three configurations: (1) only static features, (2) static features with dynamic features up to the first root LP, and (3) static features with full dynamic features extracted until the root node. As shown in Table~\ref{exp2}, models trained with only static features perform worse than using dynamic features. For example, on the NN Verification dataset, using full dynamic features improves performance up to 3.8\% compared with PD best, while using only static features results in no meaningful improvement.

\begin{mycontentbox}
    \textbf{Deeper-Stage Features Offer Superior Predictive Power} Later-stage features contain richer and more informative solver behavior signals, which are essential for improving ML-based solver configuration. Features extracted after the root node solution lead to higher predictive accuracy compared to earlier-stage features.
\end{mycontentbox}
 As shown in Table~\ref{exp2}, models trained with only static features exhibit negligible or even negative improvements, indicating the limited predictive power of structural information alone. Incorporating early-stage dynamic features (up to the first root LP) improves performance slightly, with a 0.31\% gain observed in the MIPLIB dataset. Extending to deeper-stage features (up to the root end) further amplifies this effect, achieving a 2.69\% improvement, highlighting the importance of later-stage solver behavior signals.
Figure~\ref{fig:importance} further supports this observation by illustrating the average feature importance across five datasets using the random forest. Dynamic features, particularly those captured after the root node, consistently rank higher in importance, underscoring their critical role in enhancing predictive accuracy and demonstrating that deeper-stage features provide more valuable information for ML-based solver configuration.

\begin{table}[h!]
    \centering
    \caption{Performance of Graph-based Deep Learning Methods}
    \label{exp_gnn}
    {\tiny
    \begin{tabular}{lcccccc}
        \toprule
        \textbf{Dataset} & \textbf{Method} & \textbf{Default} & \textbf{PD Best} & \textbf{Predict Time} & \textbf{Imp. Default} & \textbf{Imp. PD Best} \\
        \midrule
        \multirow{2}{*}{Independent Set}
        & GNN-MLP & 5.92 & 5.79 & 5.78 & 2.36\% & 0.17\% \\
        & L2P     & 5.92 & 5.79 & 5.73 & 3.21\% & 1.04\% \\
        \midrule
        \multirow{2}{*}{Load Balance}
        & GNN-MLP & 1123.94 & 1068.9 & 1052.2 & 6.38\% & 1.56\% \\
        & L2P     & 1123.94 & 1068.9 & 1049.8 & 6.60\% & 1.79\% \\
        \midrule
        \multirow{2}{*}{NN Verification}
        & GNN-MLP & 6.98 & 6.5 & 6.47 & 7.31\% & 0.46\% \\
        & L2P     & 6.98 & 6.5 & 6.38 & 8.60\% & 1.85\% \\
        \midrule
        \multirow{2}{*}{Set Cover}
        & GNN-MLP & 2.89 & 2.14 & 2.18 & 24.57\% & -1.87\% \\
        & L2P     & 2.89 & 2.14 & 2.15 & 25.61\% & -0.47\% \\
        \midrule
        \multirow{2}{*}{MIPLIB}
        & GNN-MLP & 46.55 & 46.55 & 46.5 & 0.11\% & 0.11\% \\
        & L2P     & 46.55 & 46.55 & 46.23 & 0.69\% & 0.69\% \\
        \bottomrule
    \end{tabular}}
\end{table}

\begin{mycontentbox}
 \textbf{Appropriate Type of Learning Method Is Important} The performance differences among different types of learning algorithms are evident. In particular, tree-based ensemble models exhibit relatively stable performance, whereas ranking-based methods demonstrate inferior results.
\end{mycontentbox}

Our analysis reveals that almost all methods can achieve positive improvements compared to the default configuration baseline. Among them, tree-based models, particularly regression models like Random Forest, deliver the best results across various datasets, consistently improving solver performance when compared to the PD best baseline. Ranker-based models perform the worst, as they focus only on relative ranking without crucial information such as solving time.

However, even with the best-performing methods, the improvements relative to the PD best baseline remain marginal, typically under 5\%. This underscores the inherent challenges in learning to optimize solver configurations effectively. The limited gains suggest that existing ML models may struggle to fully capture the complex interactions within diverse MIP instances, raising the question of whether more sophisticated learning strategies are truly necessary—or if, sometimes, all you need is the well-tuned PD best baseline.

\begin{mycontentbox}
    \textbf{Deep Learning Is Not a Panacea} Contrary to common assumptions, DL models and graph embeddings do not always outperform machine learning methods, which emphasizes that traditional models remain competitive and should not be overlooked in optimizer configuration.
\end{mycontentbox}

Random Forest and LightGBM results in  Table~\ref{exp2} consistently outperform deep learning models across multiple datasets. In some cases in Table~\ref{exp_gnn}, deep learning models even result in negative improvements. These findings demonstrate that since DL models are more sensitive to data scarcity, feature quality, and overfitting, the use of deep learning methods entails greater challenges.

\section{Open-source Toolsets BenLOC}
All components of BenLOC, including the proposed framework and evaluation criteria, are seamlessly integrated into a comprehensive open-source toolset  (https://github.com/Lhongpei/BenLOC) that supports standardized dataset management, feature extraction, model training, and performance evaluation. By providing a unified platform for ML and DL-based MIP solver configuration, BenLOC facilitates reproducible research, enables extensive benchmarking across multiple datasets, and fosters the development of advanced configuration strategies leveraging both shallow and graph-based feature representations. For more information, please refer to Appendix \ref{BenLOC}.

\section{Conclusion}
The automatic configuration of MIP optimizers remains a pivotal challenge in combinatorial optimization. In this paper, we introduce BenLOC, a standard learning-based benchmark for per-instance optimizer configuration, which ensures unbiased evaluation and robust generalization. We hope our
findings and toolsets promote future advancements in learning-based solver optimization.

\section{Acknowledgments}
The authors thank Qi Huangfu from Cardinal Operations for his insightful suggestion and assistance. This research was partially supported by the National Natural Science Foundation of China [Grants 7222500972394360, and 723943651, Natural Science Foundation of Shanghai [Grant 24ZR1421300]

\bibliographystyle{plainnat}
\bibliography{main}
\appendix
\section{Datasets}
\label{dataset}
Our approach is evaluated on different problem datasets from diverse application areas for each challenge. A problem dataset consists of a collection of MILP instances in the standard MPS file format. Each dataset was split into a training set \( D_{\text{train}} \) and a testing set \( D_{\text{test}} \), following an approximate 80-20 split. Moreover, we split the dataset by time and "optimality", which means according to the proportion of optimality for each configuration is similar in training and testing sets. This ensures a balanced representation of both temporal variations and the highest levels of configuration efficiency in our data partitions.

\subsection{Homogeneous Datasets}

The homogeneous datasets denote the datasets that comprise the same kind of problems and contain instances from only a single application. We follow the work in the literatures \cite{gasse2019exact}\cite{gasse2022machine}, and use popular datasets in our experiments. Specifically, we evaluate our approach on four homogeneous datasets:

\begin{itemize}
    \item \textbf{Set Covering: }This dataset consists of instances of the classic Set Covering Problem. Each instance requires finding the minimum number of sets that cover all elements in a universe. The problem is formulated as a MIP problem. This dataset is generated on our own and contains 3994 instances.

\item \textbf{Independent Set: }This dataset addresses the Independent Set Problem, where the goal is to find the largest set of vertices in a graph such that no two vertices are adjacent. Each instance is modeled as a MIP, with the objective of maximizing the size of the independent set. This dataset is generated on our own and contains 1604 instances.

    \item \textbf{NN Verification: }This problem is to verify whether a neural network is robust to input perturbations can be posed as a MIP \cite{cheng2017maximum}. Each input on which to verify the network gives rise to a different MIP. In this dataset, a convolutional neural network is verified on each image in the MNIST dataset, giving rise to a corresponding dataset of MIPs. This dataset contains 3104 instances.
    
    \item \textbf{Load Balancing: }This dataset is from NeurIPS 2021 Competition \cite{gasse2022machine} and it's a hard dataset. This problem deals with apportioning workloads. The apportionment is required to be robust to any worker’s failure. Each instance problem is modeled as a MILP, using a bin-packing with an apportionment formulation. This dataset contains 2286 instances.
    
\end{itemize}

The configurations used in those datasets are the same, including:

\begin{itemize}
    \item SubMipHeurLevel=0, SubMipHeurLevel=1 
    \item RootCutLevel=0, RootCutLevel=1
    \item StrongBranching=0, StrongBranching=1
    \item TreeCutLevel=0, TreeCutLevel=1
    \item DivingHeurLevel=0, DivingHeurLevel=1
    \item RoundingHeurLevel=0, RoundingHeurLevel=1
    \item MipLogLevel=2
\end{itemize} 

For each configuration, we change one configuration at a time and others are the same as the default setting of the COPT optimizer. Specifically, 'MipLogLevel=2' denotes the default configuration.

\subsection{Heterogeneous Datasets}

A heterogeneous dataset refers to a collection of data with significant variation in structure, size, and complexity, often originating from multiple application areas or problem types. This diversity poses unique challenges, as algorithms need to generalize across different characteristics and problem scales. Such datasets are particularly valuable for benchmarking solvers, as they ensure a robust evaluation of an algorithm's performance across a wide spectrum of problem instances.

In our work, we employed the heterogeneous dataset from MIPLIB 2017 \cite{gleixner2021miplib}, a well-established benchmark for evaluating MILP solvers, which is accessible at \cite{MIPLIB2017}. The dataset includes a diverse set of particularly challenging mixed-integer programming (MIP) instances, each known for its computational difficulty. These instances span various application domains, encompassing problems with different levels of complexity and solution difficulty, making MIPLIB 2017 a rigorous standard for testing the efficacy of MIP solvers in handling hard-to-solve cases.

We construct a ground dataset as follows. First, we generate ten random permutations of the elements in MIPLIB dataset, each corresponding to a seed \( s = 0, \ldots, 9 \) (where \( s = 0 \) represents the identity permutation) to expand and diversify the problem set. This expanded set is denoted as MIPLIB*, with each element representing a specific instance of the original problem set. It is important to note that, although the ten instances derived from each permutation are mathematically equivalent, they may exhibit significantly different computational behavior. For each instance, we extract information that is used to compute static features.

Next, we run the COPT solver on each instance using the configuration 'TreeCutLevel=0' and 'TreeCutLevel=1', with a time limit of 7200 seconds. The use of permutations to achieve this split is inspired by previous work \cite{berthold2022learning}.

\section{Solver Configuration}
For a specific dataset, we run the COPT solver using different configurations and collect the solver-dependent data to calculate both dynamic features and labels. Using the collected raw data, we construct a dataset of feature-label observations $D$. Suppose $P$ denotes configurations set and $F$ denotes the features, then \( D = D(P) = \{(x_{p}, y_{p}) : p \in P\} \subseteq F \times R \), which is utilized in our learning experiments.

The detailed explanation of configurations is in Table \ref{configuration}. And the possible intensity values represents: \textbf{(1) -1: }Default; \textbf{(2) 0: }Closed; \textbf{(3) 1: }Less and quick; \textbf{(4) 2: }Normal; \textbf{(5) 3: }As much as possible.

\begin{table*}[htbp]
\centering
\caption{configurations Description}
\resizebox{\textwidth}{!}{ 
\begin{tabular}{|c|c|c|c|}
\hline
\textbf{configuration Name}    & \textbf{Type}      & \textbf{Description}                              & \textbf{Values}        \\ \hline
RootCutLevel     & Integer configuration  & Strength of cuts generated at root node           & -1, 0, 1, 2, 3                  \\ \hline
TreeCutLevel      & Integer configuration  & Strength of cuts generated in the search tree     & -1, 0, 1, 2, 3                  \\ \hline
RoundingHeurLevel & Integer configuration  & Strength of Rounding heuristic                    & -1, 0, 1, 2, 3                  \\ \hline
DivingHeurLevel   & Integer configuration  & Strength of Diving heuristic                      & -1, 0, 1, 2, 3                  \\ \hline
SubMipHeurLevel   & Integer configuration  & Strength of heuristics based on sub-MIPs          & -1, 0, 1, 2, 3                  \\ \hline
StrongBranching   & Integer configuration  & Strength of Strong Branching                      & -1, 0, 1, 2, 3                  \\ \hline
\end{tabular}
}

\label{configuration}
\end{table*}
\section{Experiment Setting Details}
 \subsection{Infrastructure and Configuration}
We use 32 cores of an Intel Xeon Platinum 8469C at 2.60 GHz CPU with 512 GB RAM with 4 NVIDIA H100 GPUs. Our core functions are implemented using scikit-learn \cite{scikit-learn}, PyTorch \cite{pytorch} PyTorch Lightning and PyTorch Geometric(PyG) \cite{pyg}.

\section{Graph Based Features}
\label{gnn-feature}

\subsection{Graph Neural Networks}
Among the various types of GNN, graph attention networks (GATs) have proven to be popular and effective in practice. In our work, we adopt GATv2 \cite{gatv2}, an improved version of GAT that has been shown to be more efficient in various applications. GATv2 is an aggregation-based GNN that uses a dynamic attention mechanism to gather information from neighboring nodes connected by edges. 

The attention score \(\alpha_{ij}\) between nodes is computed as:
\begin{align*}
&\alpha_{i, j} = \frac{\exp(e_{ij})}{\sum_{k \in \mathcal{N}(i) \cup \{i\}} \exp(e_{ik})}, \text{ where } \\
&e_{ij} = 
\begin{cases} 
\mathbf{a}^{\top} \operatorname{LeakyReLU}(\boldsymbol{\Theta}_{s} \mathbf{x}_{i} + \boldsymbol{\Theta}_{t} \mathbf{x}_{j}) \text{, no edge features used} \\ 
\mathbf{a}^{\top} \operatorname{LeakyReLU}(\boldsymbol{\Theta}_{s} \mathbf{x}_{i} + \boldsymbol{\Theta}_{t} \mathbf{x}_{j} + \boldsymbol{\Theta}_{e} \mathbf{e}_{ij}) \text{, otherwise}
\end{cases}
\end{align*}

The updated feature vector for node \(i\) is then computed as:
\begin{align*}
\mathbf{x}_{i}^{\prime} = \sum_{j \in \mathcal{N}(i) \cup \{i\}} \alpha_{i,j} \boldsymbol{\Theta}_{t} \mathbf{x}_{j}
\end{align*}

Here, \(\mathcal{N}(i)\) represents the set of neighbors in the first hop of node \(i\), \(\mathbf{x}_i\) is the input feature or the embedding of node \(i\), and \(\mathbf{x}_i^\prime\) is the embedding of the output. The matrices \(\boldsymbol{\Theta}_s\), \(\boldsymbol{\Theta}_t\), and \(\boldsymbol{\Theta}_e\) are independent linear transformations applied to the features of the start node, target node, and edge features, respectively. In our work, the graph is heterogeneous with 2 types of nodes and 2 types of edge respectively, so we apply GATv2 to a heterogeneous graph by performing the aggregation process described above for each type of edge separately with independent configurations and refer to this approach as HGAT.

\subsection{Bipartite Graph Representation}

Neural networks are commonly used to extract information from structured states. In our case, the states consist of a varying number of nodes and edges, making them unsuitable for representation using fixed-size tensors. To address this, we employ graph neural networks (GNNs), which have demonstrated superior performance in the handling of structured data in several domains, such as combinatorial optimization (CO) \cite{gasse2019exact,Liu_2024} and recommendation systems \cite{gao2023surveygraphneuralnetworks}.

\section{Dynamic Features}

In this paper, we expand the dynamic features based on \cite{berthold2022learning}. In addition to the three groups of dynamic features from the Presolving, Scaling, and Global Cutting processes, we incorporated four dynamic features after solving the root LP. For some specific configurations, such as Treecutlevel, we further included seven additional dynamic features after completing the root node solution. All these features can be extracted from the COPT log file, with detailed descriptions provided in Table \ref{constraint} and \ref{static}. Moreover, the newly added features are essential, demonstrating high feature importance in the trained random forest model.

\begin{table*}[h!]
\centering
\caption{Feature and meaning of constraints and variables in the model}
\label{constraint}
\resizebox{\textwidth}{!}{ 
\begin{tabular}{|c|c|c|}
\hline
\textbf{Object} & \textbf{Feature} & \textbf{Meaning} \\
\hline
{Constraint node $\mathcal{C} = \{c_1 \dots c_m\}$} & $lb = \{lb_1 \dots lb_m\}$ & Lower-bound \\
\cline{2-3}
 & $ub = \{ub_1 \dots ub_m\}$ & Upper-bound \\
\cline{2-3}
 & $hlb = \{hlb_1 \dots hlb_m\}$ & Indicator of lower-bound \\
\cline{2-3}
 & $hub = \{hub_1 \dots hub_m\}$ & Indicator of upper-bound \\
\hline
{Variable node $\mathcal{V} = \{v_1 \dots v_n\}$} & $hlb = \{hlb_1 \dots hlb_n\}$ & Indicator of lower-bound \\
\cline{2-3}
 & $hub = \{hub_1 \dots hub_n\}$ & Indicator of upper-bound \\
\cline{2-3}
 & $c = \{c_1 \dots c_n\}$ & Coefficient in objective function \\
\cline{2-3}
 & $lb = \{lb_1 \dots lb_n\}$ & Lower-bound \\
\cline{2-3}
 & $ub = \{ub_1 \dots ub_n\}$ & Upper-bound \\
\cline{2-3}
 & $t = \{t_1 \dots t_n\}$ & Type of the variable ($\mathcal{I}, \mathcal{C}$) \\

\hline
Edge $\mathcal{E}$ & Non-zero weights in $constraints$ & \\
\hline
\end{tabular}
}

\end{table*}

\begin{table}[h!]
\centering
\caption{Feature Definitions for MIP Solver}
\label{static}
{ \small
\begin{tabular}{|c|l|}
\hline
\textbf{Feature} & \textbf{Definition} \\
\hline
\multicolumn{2}{|c|}{\textasciitilde Matrix \textasciitilde} \\
\hline
Rows & $\ln(m)$ \\
Columns & $\ln(n)$ \\
NonZeros & ratio of non-zeros, over $m \times n$ \\
Symmetries & 1 if any symmetry, 0 otherwise \\
\hline
\multicolumn{2}{|c|}{\textasciitilde Variables \textasciitilde} \\
\hline
Binaries & ratio of binary variables, over $n$ \\
Integers & ratio of integer variables, over $n$ \\
\hline
\multicolumn{2}{|c|}{\textasciitilde Constraints \textasciitilde} \\
\hline
LessThan &  \\
GreaterThan & \\
Equality &  \\
SetPartitioning &  \\
SetPacking &  \\
SetCovering & \\
Cardinality &  \\
KnapsackEquality & ratio of constraints per constraint type, over $m$ \\
Knapsack &  \\
KnapsackInteger &  \\
BinaryPacking &  \\
VariableLowerBound &  \\
VariableUpperBound & \\
MixedBinary &  \\
MixedInteger & \\
Continuous & \\
\hline
\multicolumn{2}{|c|}{\textasciitilde Scaling \textasciitilde} \\
\hline
Coefficient\_00m & $\ln(\max A'/ \min A')$ \\
RightHandSide\_00m & $\ln(\max B'/ \min B')$ \\
Objective\_00m & $\ln(\max C'/ \min C')$ \\
\hline
\multicolumn{2}{|c|}{\textasciitilde Presolving \textasciitilde} \\
\hline
PresolRows & $\ln(m)$ \\
PresolColumns & $\ln(n)$ \\
PresolIntegers & ratio of presolved integer variables, over $n$ \\
\hline
\multicolumn{2}{|c|}{\textasciitilde Global Cutting \textasciitilde} \\
\hline
DualInitialGap & $\frac{|c_d - c_l|}{\max(|c_d|, |c_l|, |c_d - c_l|)}$ \\
PrimalDualGap & $\frac{|c_p - c_l|}{\max(|c_p|, |c_l|, |c_p - c_l|)}$ \\
PrimalInitialGap & $\frac{|c_p - c_l|}{\max(|c_p|, |c_l|, |c_p - c_l|)}$ \\
GapClosed & $1 - \text{PrimalDualGap}$ \\
\hline
\end{tabular}}
\end{table}



\section{Experiment Results}

{\tiny
\begin{longtable}{lcccccccccccc}

    \caption{Performance Comparison with Different Datasets, Features and Baselines with Handcrafted Features} \label{exp2-all} \\

    \toprule
    &  \textbf{Predictor} & \textbf{Default} & \textbf{PD Best}
     & \multicolumn{3}{c}{\textbf{Static Feature}}& \multicolumn{3}{c}{\textbf{Static \& First Root Dynamic Feat.}} & \multicolumn{3}{c}{\textbf{Static \& Root End Dynamic Feat.}} \\ 
    \cmidrule(lr){5-7} \cmidrule(lr){8-10} \cmidrule(lr){11-13}
    & &   & & \textbf{Predict} & \textbf{Imp.}& \textbf{Imp.} & \textbf{Predict} & \textbf{Imp.}& \textbf{Imp.} & \textbf{Predict} & \textbf{Imp.}& \textbf{Imp.} \\
    & &  & & \textbf{Time}& \textbf{Default} & \textbf{PD Best} & \textbf{Time} & \textbf{Default} & \textbf{PD Best} & \textbf{Time} & \textbf{Default} & \textbf{PD Best} \\
    \endfirsthead

    \caption*{Continued from previous page} \\

     \toprule
    &  \textbf{Predictor} & \textbf{Default} & \textbf{PD Best}
     & \multicolumn{3}{c}{\textbf{Static Feature}}& \multicolumn{3}{c}{\textbf{Static \& First Root Dynamic Feat.}} & \multicolumn{3}{c}{\textbf{Static \& Root End Dynamic Feat.}} \\ 
    \cmidrule(lr){5-7} \cmidrule(lr){8-10} \cmidrule(lr){11-13}
    & &   & & \textbf{Predict} & \textbf{Imp.}& \textbf{Imp.} & \textbf{Predict} & \textbf{Imp.}& \textbf{Imp.} & \textbf{Predict} & \textbf{Imp.}& \textbf{Imp.} \\
    & & &  & \textbf{Time}& \textbf{Default} & \textbf{PD Best} & \textbf{Time} & \textbf{Default} & \textbf{PD Best} & \textbf{Time} & \textbf{Default} & \textbf{PD Best} \\
    \midrule
    \endhead

    \midrule
    \multicolumn{13}{r}{Continued on next page...} \\
    \endfoot

    \bottomrule
    \endlastfoot
        \midrule
        
    \multirow{13}{*}{\shortstack{Independent\\Set}} 
    &\textbf{Random Forest} & 5.92 & 5.79 &5.73& 3.21\% & 0.93\% & 5.73& 3.21\% & 0.94\% & 5.73& 3.20\% & 0.96\% \\ 
    
    & LightGBM
      & 5.92 & 5.79 &5.75&2.79\% & 0.64\% & 5.74& 2.96\% & 0.81\% & 5.73& 3.13\% & 0.98\% \\ 
    
    & XGBoost
      & 5.92 & 5.79 &5.75&2.79\% & 0.64\% & 5.74& 2.96\% & 0.81\% & 5.73& 3.13\% & 0.98\% \\
     
    & GBDT
      & 5.92 & 5.79 &5.81& 1.78\% & -0.40\% & 5.74& 2.92\% & 0.80\% & 5.58& 5.67\% & 3.58\% \\ 
    \cmidrule(lr){2-13}
    
    & LGBM Ranker  & 5.92 & 5.79 &5.86& 0.93\% & -1.26\% & 5.84& 1.27\% & -0.92\% & 5.88& 0.60\% & -1.61\% \\ 
   
    & XGB Ranker
      & 5.92 & 5.79 & 5.86 & 1.01\% & -1.21\% & 5.88 &0.55\% & -1.65\% & 5.86&1.01\% & -1.21\% \\ 
     \cmidrule(lr){2-13}
    & Rfr Clf
      & 5.92 & 5.79 & 5.85 & 1.10\% & -1.09\% & 5.69 &3.81\% & 1.67\% & 5.75&2.79\% & 0.64\% \\ 
       
    & LGBM Clf
      & 5.92 & 5.79 & 5.86 & 0.96\% & -1.23\% & 5.84 &1.35\% & -0.86\% & 5.83&1.44\%& -0.74\% \\ 

       & XGB Clf
      & 5.92 & 5.79 & 5.86 & 0.95\% & -1.24\% & 5.84 &1.27\% & -0.91\% & 5.84&1.30\%& -0.88\% \\ 

       & GBDT Clf
      & 5.92 & 5.79 & 5.81 & 1.78\% & -0.40\% & 5.80 &1.95\% & -0.23\% & 5.76&2.63\%& 0.47\% \\ 

       & KNN Clf
      & 5.92 & 5.79 & 5.79 & 2.12\% & -0.05\% & 5.77 &2.46\% & 0.29\% & 5.75&2.79\%& 0.64\% \\ 

      & SVM
      & 5.92 & 5.79 & 5.87 & 0.77\% & -1.44\% & 5.85 &1.10\% & -1.09\% & 5.77&2.46\%& 0.29\% \\ 
      
      \cmidrule[\heavyrulewidth](lr){2-13}
    &TableNet & 5.92 & 5.79 &5.75& 2.86\% & 0.72\% & 5.74& 2.93\% & 0.79\% & 5.73& 3.16\% & 0.93\% \\ 

    &MLP & 5.92 & 5.79 &5.91& 0.17\% & -2.07\% & 5.87& 0.84\% & -1.38\% & 5.84& 1.35\% & -0.86\% \\  
        
        \midrule
        \multirow{13}{*}{\shortstack{Load\\Balance}} 
         &\textbf{Random Forest}  & 1123.94 & 1068.90  & 1057.40& 5.92\% & 1.08\% & 1051.60& 6.44\% & 1.62\% & 1034.41& 8.23\% & 3.23\% \\

        & LightGBM  & 1123.94 & 1068.90  & 1059.55& 5.73\% & 0.88\% & 1057.45& 5.92\% & 1.07\% & 1030.25& 8.34\% & 3.62\% \\

        & XGBoost  & 1123.94 & 1068.90  & 1068.33& 4.95\% & 0.05\% & 1071.79& 4.64\% & -0.27\% & 1028.43& 8.50\% & 3.79\% \\

        & GBDT  & 1123.94  & 1068.90 & 1059.79& 5.71\% & 0.85\% & 1061.82& 5.53\% & 0.66\% & 1027.56& 8.58\% & 3.87\% \\ 
        \cmidrule(lr){2-13}
        
        & LGBM Ranker  & 1123.94 & 1068.90 & 1091.53&2.88\%  & -2.12\% & 1091.53&2.88\%  & -2.12\% & 1093.07& 2.75\% & -2.26\% \\ 
      
        & XGB Ranker  & 1123.94 & 1068.90  & 1084.89& 3.47\% & -1.50\% & 1085.24& 3.44\% & -1.53\% & 1086.81& 3.30\% & -1.68\% \\ 
         \cmidrule(lr){2-13}
         
        & Rfr Clf & 1123.94 & 1068.90 & 1064.70& 5.27\% & 0.39\% & 1063.92& 5.34\% & 0.47\% & 1046.43& 6.90\% & 2.10\% \\

        & LGBM Clf  & 1123.94 & 1068.90 & 1064.26& 5.31\% & 0.43\% & 1064.26& 5.31\% & 0.43\% & 1055.76& 6.07\% & 1.23\% \\ 
        
         & XGB Clf  & 1123.94 & 1068.90 & 1065.03& 5.24\% & 0.36\% & 1064.28& 5.31\% & 0.43\% & 1048.99& 6.67\% & 1.86\% \\ 
         
          & GBDT Clf  & 1123.94 & 1068.90  & 1065.03& 5.24\% & 0.36\% & 1064.07& 5.33\% & 0.45\% & 1057.45& 5.92\% & 1.07\% \\ 
          
           & KNN Clf & 1123.94 & 1068.90 & 1070.50& 4.75\% & -0.15\% & 1069.35& 4.86\% & -0.05\% & 1035.90& 7.83\% & 3.09\% \\ 

            & SVM  & 1123.94 & 1068.90  & 1064.38& 5.30\% & 0.42\% & 1064.30& 5.31\% & 0.43\% & 1064.07& 5.33\% & 0.45\% \\ 
            
       \cmidrule[\heavyrulewidth](lr){2-13}
       &TableNet  & 1123.94 & 1068.90  & 1064.70& 5.27\% & 0.39\% & 1048.99& 6.68\% & 1.86\% & 1035.90& 7.84\% & 3.09\% \\ 
       &MLP & 1123.94 & 1068.90  & 1086.70& 3.31\% & -1.67\% & 1082.23& 3.71\% & -1.25\% & 1075.49& 4.31\% & -0.62\% \\ 
        
        \midrule
        
        \multirow{13}{*}{\shortstack{NN\\Verification}}
          &\textbf{Random Forest}  & 6.98 & 6.50    & 6.59& 5.59\% & -1.58\% & 6.58& 5.73\% & -1.40\% & 6.26& 10.32\% & 3.69\% \\

        & LightGBM  & 6.98 & 6.50    & 6.62& 5.13\% & -1.89\% & 6.63& 4.99\% & -2.05\% & 6.25& 10.43\% & 3.80\% \\

        & XGBoost  & 6.98 & 6.50    & 6.62& 5.13\% & -1.89\% & 6.64& 4.84\% & -2.20\% & 6.29& 9.86\% & 3.19\% \\

       & GBDT  & 6.98 & 6.50    & 6.67& 4.41\% & -2.66\% & 6.63& 4.99\% & -2.05\% & 6.28& 10.00\% & 3.34\% \\ 
        \cmidrule(lr){2-13}
        
        & LGBM Ranker  & 6.98 & 6.50    & 6.72& 3.70\% & -3.43\% & 6.72& 3.70\% & -3.43\% & 6.74& 3.41\% & -3.74\% \\

        & XGB Ranker  & 6.98 & 6.50    & 6.71& 3.84\% & -3.28\% & 6.74& 3.41\% & -3.74\% & 6.70& 3.98\% & -3.12\% \\ \cmidrule(lr){2-13}
       
        & Rfr Clf  & 6.98 & 6.50    & 6.67& 4.41\% & -2.66\% & 6.61& 5.27\% & -1.74\% & 6.37& 8.71\% & 1.95\% \\

        & LGBM Clf & 6.98 & 6.50    & 6.69& 4.13\% & -2.97\% & 6.60& 5.42\% & -1.59\% & 6.43& 7.85\% & 1.03\% \\ 
        
         & XGB Clf & 6.98 & 6.50    & 6.66& 4.56\% & -2.51\% & 6.59& 5.56\% & -1.43\% & 6.40& 8.28\% & 1.49\% \\ 
         
          & GBDT Clf & 6.98 & 6.50    & 6.64& 4.84\% & -2.20\% & 6.58& 5.70\% & -1.28\% & 6.41& 8.14\% & 1.34\% \\ 
          
           & KNN Clf & 6.98 & 6.50    & 6.65& 4.70\% & -2.36\% & 6.58& 5.70\% & -1.28\% & 6.38& 8.57\% & 1.80\% \\ 

            & SVM & 6.98 & 6.50    & 6.73& 3.55\% & -3.59\% & 6.67& 4.41\% & -2.66\% & 6.63& 4.99\% & -2.05\% \\ 
            
        \cmidrule[\heavyrulewidth](lr){2-13}

        &TableNet  & 6.98 & 6.50    & 6.62& 5.07\% & -1.95\% & 6.59& 5.59\% & -1.58\% & 6.28& 10.12\% & 3.39\% \\ 

        &MLP  & 6.98 & 6.50    & 6.78& 2.87\% & -4.31\% & 6.62& 5.15\% & -1.85\% & 6.56& 6.02\% & -0.92\% \\
        
         \midrule
         
        \multirow{13}{*}{Set Cover} 
        &Random Forest & 2.89 & 2.14 &2.14&26.19\% & -0.03\% & 2.14& 26.20\% & -0.02\% & 2.13& 26.30\% & 0.04\% \\

        & LightGBM
          & 2.89 & 2.14  &2.14&26.19\% & -0.03\% & 2.14& 26.20\% & -0.02\% & 2.15& 26.05\% & -0.12\% \\

        & XGBoost
          & 2.89 & 2.14 &2.14&26.23\% & -0.01\% & 2.14& 26.21\% & -0.03\% & 2.14& 26.21\% & -0.04\% \\

        & \textbf{GBDT}
          & 2.89 & 2.14  &2.14&26.19\% & -0.03\% & 2.14& 26.21\% & -0.02\% & 2.13& 26.44\% & 0.30\% \\
          
        \cmidrule(lr){2-13}
        & LGBM Ranker  & 2.89 & 2.14 &2.58& 10.75\% & -20.97\% & 2.58& 10.71\% & -21.02\% & 2.58& 10.75\% & -20.97\% \\

        & XGB Ranker
          & 2.89 & 2.14  & 2.59 & 10.68\% & -21.06\% & 2.58 &10.75\% & -20.96\% & 2.58 &10.61\% & -21.17\% \\ 
          \cmidrule(lr){2-13}
        
        & Rfr Clf
          & 2.89 & 2.14  &2.14&26.23\% & -0.00\% & 2.14& 26.21\% & -0.03\% & 2.14& 26.15\% & -0.11\% \\

        & LGBM Clf
          & 2.89 & 2.14  &2.14&26.19\% & -0.03\% & 2.13& 26.22\% & -0.01\% & 2.13& 26.24\% & -0.00\% \\

         & XGB Clf
          & 2.89 & 2.14  & 2.14 & 26.18\% & -0.07\% & 2.14 & 26.22\% & -0.01\% & 2.14 & 26.23\% & -0.00\% \\

          & GBDT Clf
          & 2.89 & 2.14  & 2.14 & 26.19\% & -0.06\% & 2.14 & 26.22\% & -0.01\% & 2.14 & 26.23\% & -0.00\% \\

          & KNN Clf
          & 2.89 & 2.14  & 2.14 & 26.02\% & -0.29\% & 2.14 & 26.04\% & -0.26\% & 2.14 & 26.19\% & -0.05\% \\

           & SVM
          & 2.89 & 2.14  & 2.14 & 26.21\% & -0.03\% & 2.13 & 26.25\% & 0.00\% & 2.13 & 26.25\% & 0.00\% \\
          
          \cmidrule[\heavyrulewidth](lr){2-13}
          
          &TableNet & 2.89 & 2.14  &2.14 &26.18\% & -0.07\% & 2.14& 26.19\% & -0.06\% & 2.13& 26.28\% & 0.02\% \\

          &MLP & 2.89 & 2.14  &2.36 &18.34\% & -10.28\% & 2.27& 21.45\% & -6.75\% & 2.19 & 24.22\% & -2.34\% \\
          
           \midrule
        \multirow{13}{*}{MIPLIB}
         &\textbf{Random Forest}  & 46.55 & 46.55 & 46.62 & -0.15\% & -0.15\% & 46.50 & 0.12\% & 0.12\% & 45.39 & 2.49\% & 2.49\%\\
         
         &LightGBM  & 46.55 & 46.55 & 47.55 & -2.31\% & -2.31\%  & 46.41 & 0.31\% & 0.31\% & 45.32 & 2.69\% & 2.69\%\\

         &XGBoost  & 46.55 & 46.55 & 47.30 & -1.65\% & -1.65\% & 46.52 & 0.06\% & 0.06\% & 46.32 & 0.41\% & 0.41\%\\

         &GBDT  & 46.55 & 46.55 & 47.08 & -1.14\% & -1.14\% & 46.48 & 0.15\% & 0.15\% & 45.96 & 1.27\% & 1.27\%\\
         
\cmidrule(lr){2-13}
         &LGBM Ranker  & 46.55 & 46.55 & 47.62 & -2.30\% & -2.30\% & 47.50 & -2.05\% & -2.05\% & 46.88 & -0.71\% & -0.71\%\\

         &XGB Ranker  & 46.55 & 46.55 & 47.64 & -2.34\% & -2.34\% & 47.23 & -1.46\% & -1.46\% & 46.97 & -0.90\% & -0.90\%\\
         \cmidrule(lr){2-13}

         &Rfr Clf  & 46.55 & 46.55 & 48.10 & -3.43\% & -3.43\% & 46.49 & 0.13\% & 0.13\% & 46.05 & 0.96\% & 0.96\%\\

         &LGBM Clf  & 46.55 & 46.55 & 47.80 & -2.69\% & -2.69\% & 46.98 & -0.99\% & -0.99\% & 46.08 & 0.90\% & 0.90\%\\

         & XGB Clf & 46.55 & 46.55  & 47.92 & -2.94\% & -2.94\% & 46.69 & -0.34\% & -0.34\% & 46.20 & 0.75\% & 0.75\% \\
         
 & GBDT Clf & 46.55 & 46.55  & 47.60 & -2.27\% & -2.27\% & 46.64 & -0.19\% & -0.19\% & 46.10 & 0.97\% & 0.97\% \\
 
 & KNN Clf & 46.55 & 46.55  & 47.83 & -2.74\% & -2.74\% & 46.50 & 0.11\% & 0.11\% & 46.15 & 0.86\% & 0.86\% \\

& SVM & 46.55 & 46.55 & 47.29 & -1.59\% & -1.59\% & 46.65 & -0.21\% & -0.21\% & 46.24 & 0.67\% & 0.67\% \\

    \cmidrule[\heavyrulewidth](lr){2-13}
    & TableNet  & 46.55 & 46.55 & 46.61 & -0.13\% & -0.13\% & 46.52 & 0.06\% & 0.06\% & 45.42 & 2.43\% & 2.43\%\\    

     & MLP  & 46.55 & 46.55 & 48.47 & -4.13\% & -4.13\% & 47.00 & -0.96\% & -0.96\% & 46.78 & -0.49\% & -0.49\%\\    

\end{longtable}}

\section{Open source Toolsets BenLOC}
\label{BenLOC}
\subsection{Open Datasets}
Our library ships with most frequently using benchmark collections referenced in this paper, including \texttt{setcover}, \texttt{indset}, \texttt{nn\_verification}, \texttt{load\_balancing} and \texttt{MIPLIB}, refer to Appendix \ref{dataset} for more details.

\subsection{Key Features}

\par\textbf{Flexible API \& End-to-End Pipeline:} Our framework delivers a fully integrated, end-to-end pipeline for machine-learning–driven solver configuration. It conforms to the familiar \texttt{fit}/\texttt{predict} interface, allowing practitioners to plug in any regression or classification model-ranging from classical algorithms such as \texttt{RandomForestRegressor} and \texttt{LinearRegression} to custom estimators defined by the user-without touching the surrounding workflow. For deep-learning use cases, first-class PyTorch support is provided: datasets are wrapped as \texttt{torch.utils.data.Dataset} objects, models can be defined as subclasses of \texttt{torch.nn.Module}, and training loops-including checkpointing, learning-rate scheduling, and GPU management-are orchestrated automatically.

\par\textbf{Integrated MILP Datasets \& Preprocessing:} To facilitate rapid experimentation on Mixed-Integer Linear Programs (MILPs), the framework includes a curated collection of benchmark datasets, notably \texttt{setcover}, \texttt{indset} (maximum independent set), and \texttt{nn\_verification}. Each instance is parsed and transformed into rich feature matrices capturing characteristics such as constraint density, coefficient distributions, and variable–constraint incidence. Users may further augment these representations with proprietary features in collaboration with solver developers-for example, presolve elimination rates, cut-generation frequencies, or solver-specific diagnostics-to capture domain-specific insights.

\par\textbf{Bipartite Graph Representation Learning:} Recognizing the inherent combinatorial structure of MILP instances, our system incorporates representation-learning modules based on a bipartite-graph view: variables and constraints form two disjoint node sets, with edges annotated by coefficient magnitudes. Through configurable message-passing layers, the framework learns embeddings that summarize global instance features, which are then fed into downstream prediction heads for runtime estimation, configuration selection, or cut-strength prediction. Empirical studies demonstrate that this graph-based approach often outperforms flat-vector feature sets, particularly when modeling higher-order interactions.

\subsection{Core Functionalities}

\par\textbf{configuration Configuration:} The \texttt{Params} class centralizes configuration of label scaling, baseline selection, and geometric-mean normalization for consistent performance evaluation. Users specify target metrics, weightings, and baseline solvers in a single, coherent interface.

\par\textbf{Flexible Training:} Model training is fully customizable: all scikit‑learn estimators are supported out of the box, and a simple adapter pattern allows seamless integration of user‑defined models or PyTorch networks. The same \texttt{fit}/\texttt{predict} contract governs both classical and deep models, ensuring consistency across experiments.

\par\textbf{Dataset Management:} Dataset management is handled transparently via \texttt{pandas.DataFrame} objects. Preloaded benchmark collections are available immediately, while user‑supplied data can be ingested, shuffled, stratified, and split into training, validation, and test folds with a single API call. Built‑in preprocessing routines extract raw features, normalize distributions, and handle missing data automatically.

\par\textbf{Performance Evaluation:} Performance evaluation leverages the \texttt{evaluate()} method, which returns detailed comparison tables against baselines and Per-Instance best configurations. Metrics include absolute speedups, relative improvements, and statistical significance tests where applicable, giving users clear insight into model impact.

\par\textbf{Hyperconfiguration Tuning:} Finally, hyperconfiguration tuning and model selection are supported via both cross‑validation and Bayesian optimization. Predefined search spaces are provided for Random Forests, LightGBM, and Gradient Boosted Decision Trees (GBDT), while users may define custom configuration ranges or objective functions. Automated early stopping and parallel evaluation make large‑scale searches efficient and reproducible.
\subsubsection{Example Usage}
An example usage of the BenLOC framework is provided below:
\begin{verbatim}
# Step 1: Initialize the model with configurations
from BenLOC.params import Params
params = Params(default="MipLogLevel-2", label_type="log_scaled")
model = BenLOC(params)

# Step 2: Set a machine learning model (e.g., Random Forest)
from sklearn.ensemble import RandomForestRegressor
model.set_trainer(RandomForestRegressor(verbose=1))

# Step 3: Load and preprocess a dataset
model.load_dataset("setcover", processed=True, verbose=True)

# Step 4: Split the dataset into training and testing sets
model.train_test_split(test_size=0.2)

# Step 5: Train the model
model.fit()

# Step 6: Evaluate the model
evaluation_results = model.evaluate()
print(evaluation_results)
\end{verbatim}


\end{document}